\documentclass[final]{template}

\usepackage{enumitem}
\setlist{topsep=0.1em, itemsep=0.1em}

\usepackage{tocbasic}
\DeclareTOCStyleEntry[
beforeskip=.2em plus 1pt,% default is 1em plus 1pt
pagenumberformat=\textbf
]{tocline}{section}

\renewcommand{\theenumi}{\Alph{enumi}}

\begin{document}

\setcounter{page}{1}

\date{}

\title{\vspace{-2cm}Choices of HKR isomorphisms}
\author{Marco Robalo\footnote{Sorbonne Université,  Institut de Mathématiques de Jussieu - Paris Rive Gauche, marco.robalo@imj-prg.fr}  }

\normalem
\vspace{-1cm}
\maketitle
\vspace{-1cm}
\begin{abstract}
We record the fact that the set of chain-level multiplicative HKR natural equivalences defined simultaneously for all derived schemes, functorialy splitting the HKR-filtration and rendering the circle action compatible with the de Rham differential,  is, via Cartier duality, in a natural bijection with the set of filtered formal exponential maps $ \widehat{\mathbb{G}_a}\to \widehat{\mathbb{G}_m}$. In particular, when the base $\field$ is a field of characteristic zero, the set of choices is $\field^\ast$.
\end{abstract}

\tableofcontents

\section{Introduction}
  Let $\field$ be a commutative ring. The Hochschild-Kostant-Rosenberg (HKR) theorem \cite{MR0142598} establishes for any smooth $\field$-scheme $X=\Spec(R)$ an identification of the Hochschild homology groups $\HH_i(R/\field):=\Tor^i_{R\otimes_\field R}(R,R)$  with the modules of $i$-differential forms  $\Omega_{X/\field}^i$, given by the anti-symmetrization map
  \[
\Omega^i_{R/\field}\to \HH_i(R/\field),\,\,\,\,\,\,\, r_0.dr_1\wedge \cdots \wedge dr_i\mapsto \sum_{\sigma\in \Sigma_i} (-1)^{\sign(\sigma)} [r_0\otimes r_{\sigma(1)}\otimes\cdots \otimes r_{\sigma(i)}].
  \]
 The groups $\mathrm{HH}_i(R/\field)$ are actually defined for every derived $\field$-algebra $R\in \SCR_\field$\footnote{Here $\SCR_\field$ denotes the \icategory of simplicial commutative rings. See \cite[\S 25.1.1]{lurie-sag}} as the homology groups of the derived tensor product of $\field$-algebras  $ \HH(R/\field):=R\derivedtensor{ R\derivedtensor{\field}R} R$ where $R$ is seen as an $R\derivedtensor{\field}R$-algebra using the multiplication map $R\derivedtensor{\field} R\to R$. In particular, this shows that $\HH(R/\field)$ carries the structure of an object in $\SCR_\field$. Also, for a general $R\in \SCR_\field$ we replace $\Omega^1_{R/\field}$ by the \emph{cotangent complex} $\cotangent_{R/\field}$ and independently of the characteristic of $\field$, we have the HKR filtration  on $\HH(R/\field) $ that has $(\Lambda^i \cotangent_{R/\field})[i]$ as 
   associated graded piece of weight $i$ (see \cite[IV. 4.1]{Nikolaus2018}).
 When $\field$ is a field with $\characteristic(\field)=0$, the anti-symmetrization map  induces a splitting of the HKR filtration and gives a $\field$-linear quasi-isomorphism \cite[Prop. 1.3.16, Remark 3.2.3, Prop. 5.4.6]{MR1217970}
	\begin{equation}
	\label{quasi-iso-hkr-formula}\begin{tikzcd}\HH(R/\field)\ar{r}[swap]{\sim}& \bigoplus_{i=1}^n \, (\Lambda^i \cotangent_{R/\field})[i]\end{tikzcd}
\end{equation} 
\noindent  Derived geometry \cite{MR2862069, MR2928082} offers another perspective:  since in  $\SCR_\field$, derived tensor products are pushouts, we find an equivalence in $\SCR_\field$, $\HH(R/k)\simeq R\derivedtensor{\field} \circle$ that presents $\HH(R/\field)$ as the derived ring of functions $\structuresheaf_{\loops X}$ on the derived loop scheme $\loops X:=\derivedMap(\circle, X)$ where by $\derivedMap$ we mean the derived mapping scheme cf. \cite[\S 3.2]{MR3285853}.\\

 \noindent Similarly, derived geometry offers a geometric incarnation for $\bigoplus_{i=1}^n \, (\Lambda^i \cotangent_{R/\field})[i]$ as the derived ring of functions of the shifted tangent bundle  $\tangentstack[-1]X$=$\Spec(\Sym^{\Delta}(\cotangent_{X/\field}[1]))$ where $\Sym^\Delta$ corresponds to the symmetric algebra construction in the setting of commutative simplicial rings, with the $\Gm{\field}$-action corresponding to the natural grading.  When $\field$ is a $\rationals$-algebra, for any affine derived scheme $X$, the results of \cite{MR2862069, MR2928082}  provide an isomorphism of derived schemes functorial in $X$

\begin{equation}
	\label{isoderivedschemesHKR}
\xymatrix{\tangentstack[-1]X\ar[r]^-{\sim}&\loops X}
\end{equation}
\noindent that recovers a quasi-isomorphism of the type \eqref{quasi-iso-hkr-formula} after passing to global functions. However, it is unclear if the equivalence obtained through derived geometry coincides with the anti-symmetrization map of \eqref{quasi-iso-hkr-formula}.

\begin{observation}\label{HKRwithconnection}
When $\characteristic(\field)=0$, Kapranov \cite{Kapranov1999} explains another way to produce HKR isomorphisms \eqref{quasi-iso-hkr-formula} by considering smooth schemes $X$ with a torsion-free flat connection $\nabla$ on their tangent bundle. In this case the connection provides a formal exponential  $\exp^\nabla: \widehat{\tangentstack X}\simeq  \widehat{\Delta}$ establishing an isomorphism between the formal completion of $X\times X$ along the diagonal and the formal completion of $\tangentstack X$ along the zero section. Passing to the self-intersections with $X$, we obtain another equivalence of derived schemes of the type  \eqref{isoderivedschemesHKR}.
\end{observation}

\section{The space of functorial HKR isomorphisms}

\noindent  The goal of this note is to prove  \cref{maincorollarynote} below, describing the collection of  HKR isomorphisms \eqref{quasi-iso-hkr-formula}.  We start by noticing though, that without further assumptions, this space can be significantly large:  \cref{HKRwithconnection} shows that every torsion-free connection on a scheme $X$ provides one, and the space of connections is affine. But clearly, connection-induced HKR isomorphisms are not functorial unless the maps preserve the connection.  Therefore we will only consider chain level HKR equivalences enhanced with:

\begin{enumerate}[itemsep=0.01cm]

	\item  functoriality for all derived  $\field$-rings as part of a natural equivalence  of \ifunctors  on the \icategory of derived schemes 

	\begin{equation}\label{HKRnaturaltransformation}\xymatrix{\tangentstack[-1](-)\ar[r]^{\sim}&\loops (-)}\end{equation}

	\item    functorial splittings of the HKR-filtration;
	\item functorial  matchings of the circle action on loop spaces with the de Rham differential on forms .

\end{enumerate}

\begin{observation}
	In particular, chain level HKR-equivalences as in (i) are automatically multiplicative by passing to the derived rings of functions in \eqref{HKRnaturaltransformation}.
\end{observation}

\noindent Before formulating our main result we must first describe how derived  geometry helps combining the structures in (i)-(iii), culminating with  \cref{spaceHKRisomorphisms} below.  Using the formalism of affine stacks \cite{MR2244263}, it is shown in the combined results of \cite{MR4444269,Toen2020c, Moulinos2021} that over any commutative ring $\field$ there exists a  flat affine filtered abelian group stack (underived)
\[\Filcircle\to \Filstack{\field}\]
\noindent which we call the filtered circle, and such that for any derived scheme $X$, the relative derived mapping stack
\[
\derivedMap_{\Filstack{\field}}(\Filcircle, X\times \Filstack{\field})\to \Filstack{\field}
\]
\noindent provides the HKR-filtration on the derived loop space $\loops X$, with associated graded stack given by $(\tangentstack[-1]X)/\Gm{\field}\to \B\Gm{\field}$.  More precisely, it describes the HKR-filtration with its multiplicative structure as the universal filtered algebra with an action of the filtered circle. As a consequence, asking for HKR-isomorphisms (i)-(iii), is to ask for \emph{splittings} of the filtered circle compatible with the group structure. Let us then recall the construction of a split filtered stack associated to a graded stack:

\begin{construction}
	\label{splitfiltstackconstruction}
Let $q:\Filstack{\field}\to \B\Gm{\field}$ be the map induced by the projection $\affine^1_\field \to \Spec(\field)$ and let $Y$ be a stack endowed with a $\Gm{}$-action. Take $Z=[Y/\Gm{}]\to \B\Gm{\field}$. We define the associated \emph{split} filtered stack $Z^{\split}\to \Filstack{\field}$ to be the pullback
	\[
	\begin{tikzcd}[column sep=small, row sep=small]
		Z^{\split}\ar{d}\ar{r}\tikzcart& Z\ar{d}\\
		\Filstack{\field}\ar{r}& \B\Gm{\field}
	\end{tikzcd}
	\]
	\noindent By construction, it is equivalent to the quotient stack $[Y\times \affineline{}/\Gm{}]$ where we let $\Gm{}$ act on the product coordinate-wise. The associated graded stack $(Z^\split)^\gr$ is canonically equivalent to $Z$ because $q$ is a right inverse to the inclusion $0:\B\Gm{\field}\to \Filstack{\field}$. Finally, when $S\to \Filstack{\field}$ is a filtered stack, we denote by $S^\triv:=(S^{\gr})^{\split}$ the associated split filtered stack where $S^\gr$ is the pullback of $S$ along the inclusion $\B\Gm{\field}\to \Filstack{\field}$.

\end{construction}

\begin{observation}
\noindent Since \cref{splitfiltstackconstruction} is monoidal with respect to cartesian products, $(\Filcircle)^\triv$ is still a group object. 
\end{observation}

\noindent We can finally formulate how to combine the enhanced structures of (i)-(iii) as part of a single object:

\begin{definition}\label{spaceHKRisomorphisms}
We define the set of  chain-level HKR-isomorphisms enhanced with (i)-(iii) as the set of connected components of the mapping space of invertible maps of group (higher) stacks
\[
\Map_{\mathrm{group}, \Filstack{\field}}^{\mathrm{inv}}\left(\Filcircle \,,\, (\Filcircle)^\triv\right)
\]
\noindent ie, universal splittings of the HKR filtration compatible with the action of the filtered circle.
\end{definition}

\begin{observation} Given a splitting $\Filcircle\simeq (\Filcircle)^\triv$ as in \cref{spaceHKRisomorphisms} we obtain the associated HKR-natural transformation \eqref{HKRnaturaltransformation} by pre-composition with the relative derived mapping spaces over $\Filstack{\field}$
\[
\begin{tikzcd}\derivedMap_{\Filstack{\field}}((\Filcircle)^\triv, X\times \Filstack{\field})\ar{r}{\sim}& \derivedMap_{\Filstack{\field}}(\Filcircle, X\times \Filstack{\field})
	\end{tikzcd}
\]
\noindent and extracting the fibers over $1:\Spec(\field)=[\Gm{\field}/\Gm{\field}]\to \Filstack{\field}$.
\end{observation}

\noindent We  state our main result:

\begin{theorem}\label{maincorollarynote}
	Let $\field$ be a field. Then, the set of chain-level HKR equivalences enhanced with (i)-(iii) (cf. \cref{spaceHKRisomorphisms}) is in bijection with the set of formal exponentials, ie, group homomorphisms of formal groups, 
	\[\Hom_{\fgps}(\formalGa{\field}, \formalGm{\field})\simeq \begin{cases}\field^\ast & \text{ if  } \characteristic(\field)=0,\\ \emptyset& \text{otherwise}.\end{cases}\]
\end{theorem}
\section{Proof of \cref{maincorollarynote}}

\noindent We are interested in computing $\pi_0$ of the space in \cref{spaceHKRisomorphisms}.  Thanks to  \cite[Theorem 1.8]{Moulinos2021} we have an explicit formula for the filtered group circle in terms of relative Cartier duality over $\Filstack{\field }$
\[
\Filcircle\simeq \B_{\Filstack{\field}}( \Deformationspace^\vee)
\]
\noindent where $\Deformationspace\to \Filstack{\field}$ is the formal group scheme over $\Filstack{\field }$ given by the total space of the deformation to the normal bundle at the unit from the formal group $\formalGm{\field}$ to its lie algebra $\formalGa{\field}$ (cf. \cite[Construction 5.6, Proposition 5.12, Theorem 1.6]{Moulinos2021}.). Here, relative Cartier duality is given by the $\Filstack{\field}$-relative construction of \cite[37.3.4]{MR2987372}:
\[(-)^\vee:= \Hom_{\fgps}(- , \formalGm{\field})\]
\noindent (the hom is taken inside the category of classical formal group schemes, not as derived schemes) and $\formalGm{\field}$ is the multiplicative formal group. Since the construction of Cartier duality is the relative one, we can freely interchange 
\[
(\Deformationspace^\vee)^\triv \simeq (\Deformationspace^\triv)^\vee, \hspace{1cm}\B_{\Filstack{\field}}(\Deformationspace^\vee)^\triv \simeq \B_{\Filstack{\field}}((\Deformationspace^\vee)^\triv)
\]
\noindent As a consequence, the space of HKR-isomorphisms of \cref{spaceHKRisomorphisms} is equivalent to
\[\Map_{\mathrm{group}, \Filstack{\field}}^{\mathrm{inv}}\left( \B_{\Filstack{\field}}( \Deformationspace^{\vee})\,,\,  \B_{\Filstack{\field}} ((\Deformationspace^\triv)^{\vee})\right)\]
\noindent Since all group stacks being used are abelian, the Eckmann–Hilton delooping at the unit provides a map 
\[
\begin{tikzcd}[row sep=small, column sep=small]\Map_{\mathrm{group}, \Filstack{\field}}^{\mathrm{inv}}\left( \B_{\Filstack{\field}}( \Deformationspace^{\vee})\,,\,  \B_{\Filstack{\field}} ((\Deformationspace^\triv)^{\vee})\right)\ar{d}{\Omega_\ast}\\\Map_{\mathrm{group}, \Filstack{\field}}^{\mathrm{inv}}\left(  \Deformationspace^{\vee}\,,\,  (\Deformationspace^\triv)^{\vee}\right)\end{tikzcd}
\]
\noindent which induces an isomorphism of $\pi_0$ with inverse given by the $\B$-construction. \\

\noindent Finally, we consider the map  induced by the functor of Cartier duality 
\begin{equation}\label{ff} \begin{tikzcd}[row sep=small]\Map_{\mathrm{group}, \Filstack{\field}}^{\mathrm{inv}}\left(  \Deformationspace^{\vee}\,,\,  (\Deformationspace^\triv)^{\vee}\right)\\
\Map_{\fgps, \Filstack{\field}}^{\mathrm{inv}}\left( \Deformationspace^\triv, \Deformationspace\right)\ar{u}{(-)^\vee}\end{tikzcd}\end{equation}
\noindent which is an equivalence, thanks to the fully faithfulness of Cartier duality \cite[Const 3.7, Prop 3.12, Const 3.17, Prop. 3.19]{Moulinos2021}. Here, 
$\fgps$ denotes the category of relative smooth formal groups. Notice that, independently of $\characteristic(\field)$, both mapping spaces in \eqref{ff} are discrete. Moreover, thanks to \cite[1.4.2 and 1.4.5]{MR3665552} we can either see the last mapping space as maps of prestacks or as continuous maps.\\

\noindent Since $\Deformationspace\to \Filstack{\field}$ is a smooth formal group relative to $\Filstack{}$, we can identify the trivial filtration $\Deformationspace^{\triv}\to \Filstack{\field}$ with the affine linear formal group associated to its relative Lie algebra. In particular, following \cref{splitfiltstackconstruction}, it is given by the constant family over $\Filstack{\field}$
\[
\Deformationspace^{\triv}\simeq [(\formalGa{\field}\times \affineline{\field})/\Gm{\field}]
\]
\noindent In conclusion, we have shown that the set of functorial HKR-isomorphisms as in \cref{spaceHKRisomorphisms} is in bijection with the set of filtered formal exponentials
\[
\Map_{\fgps, \Filstack{\field}}^{\mathrm{inv}}\left( [(\formalGa{\field}\times \affineline{\field})/\Gm{\field}]\,,\, \Deformationspace\right)
\]
\begin{observation}
	By extracting the underlying groups of the filtration (ie, the fibers over $1$ in $\Filstack{\field}$) we find a map
	
	\begin{equation}\label{heightreasons}
	\Map_{\fgps, \Filstack{\field}}^{\mathrm{inv}}\left( [(\formalGa{\field}\times \affineline{\field})/\Gm{\field}]\,,\, \Deformationspace\right)\to \Map_{\fgps }^{\mathrm{inv}}\left(\formalGa{\field}, \formalGm{\field}\right)
	\end{equation}
	
\noindent 	By height reasons, since $\formalGa{\field}$ is of height $\infty$ and $\formalGm{\field}$ is of height 1, the target of \eqref{heightreasons} is empty when $\field$ is of $\characteristic(p)>0$. Therefore, so is the source of \eqref{heightreasons}.
	
\end{observation}

\noindent Finally, assume $\characteristic(\field)=0$.  The relative exponential map (see for instance \cite[Exposé VIIB - \S3]{demazure1970schemas} or \cite[Chapter 7, Cor. 3.2.2]{MR3701353})
  defines an isomorphism of filtered formal group schemes
\[\begin{tikzcd}[column sep=small][(\formalGa{\field}\times \affineline{\field})/\Gm{\field}]\ar{r}{\exp_{\mathrm{rel}}}[swap]{\sim}& \Deformationspace
\end{tikzcd}\]
\noindent Composition with $\exp_{\mathrm{rel}}$ defines a bijection

\begin{equation}\label{equationlastexponential}\begin{tikzcd}[row sep=small] \Map_{\fgps, \Filstack{\field}}^{\mathrm{inv}}\left( [(\formalGa{\field}\times \affineline{\field})/\Gm{\field}]\,,\, \Deformationspace\right)\\
\Map_{\fgps, \Filstack{\field}}^{\mathrm{inv}}\left( [(\formalGa{\field}\times \affineline{\field})/\Gm{\field}]\,,\, [(\formalGa{\field}\times \affineline{\field})/\Gm{\field}]\right)
\ar{u}[anchor=south, rotate=90]{\sim}[swap]{\,\,\exp_{\mathrm{rel}}\,\circ\, -}
\end{tikzcd}
\end{equation}
\noindent Let us compute the last space: since $\characteristic(\field)=0$, the category of formal groups relative to $\Filstack{\field}$ is equivalent to the category of Lie algebra objects in $\Qcoh(\Filstack{\field})$ \cite[Chapter 7]{MR3701353}. The Lie algebra associated to  $[(\formalGa{\field}\times \affineline{\field})/\Gm{\field}]$ is the structure sheaf $\structuresheaf_{\Filstack{\field}}(1)$ with the weight-$(1)$ action of $\Gm{\field}$, endowed with the abelian Lie bracket (see \cite[\S 5]{Moulinos2021}).  Since $\Qcoh(\Filstack{\field})$ is symmetric monoidal equivalent to filtered $\field$-modules $\Fil(\Mod_\field)$ \cite{Tasos}, $\structuresheaf_{\Filstack{\field}}(1)$ corresponds to the abelian Lie algebra given by $\field(1)$.  It follows that
\[\begin{tikzcd}[row sep=small]
	\Map_{\fgps, \Filstack{\field}}^{\mathrm{inv}}\left( [(\formalGa{\field}\times \affineline{\field})/\Gm{\field}]\,,\, [(\formalGa{\field}\times \affineline{\field})/\Gm{\field}]\right)\\ \ar[u, "\sim\,\," labl, ]\pi_0\, \Map_{\mathrm{Lie}, \Fil(\Mod_\field)}^{\mathrm{inv}}( \field(1), \field(1))=\field^\ast
	\end{tikzcd}
\]
\noindent In particular, the map \eqref{equationlastexponential} sends  $\lambda\in \field^\ast$ to $\exp(\lambda. (-))$, thus proving \cref{maincorollarynote}.

\begin{remark}
	\cref{maincorollarynote} describes the space of \emph{group splittings} of the filtered circle as exponentials (cf. \cref{spaceHKRisomorphisms}). The results of \cite{Moulinos2022} show that even in characteristic zero, the filtered circle does not admit splittings as a \emph{pointed cogroup} 
	with co-multiplication given by the pinch map. The universal obstruction is the Todd class. Recall  that the splitting principle for algebraic $\Ktheory$-theory implies that the collection of Chern characters from $\Ktheory$-theory to de Rham cohomology coincides with the collection of exponential maps - see \cite[Lemma 5.5]{MR3338682}. In summary, the existence of  group splittings of $\Filcircle$ allows the Chern characters to exist, and the fact that none of those are cogroup splittings, imposes the Grothendieck-Riemann-Roch theorem.
	\end{remark}

\noindent \textbf{Acknowledgements:}
I thank Mauro Porta, Bertrand To\"en, Tasos Moulinos and Nick Rozenblyum for discussions about this short note.

{\tiny
\printbibliography}

\end{document}